# A Generalization of Stationary AR(1) Schemes

**Satheesh, S**
Neelolpalam, S. N. Park Road
Trichur – 680 004, India.
*ssatheesh1963@yahoo.co.in*

**Sandhya, E**
Department of Statistics, Prajyoti Niketan College
Pudukkad, Trichur – 680 301, India**.**
*esandhya@hotmail.com*

and

**Sherly, S**
Department of Statistics, Vimala College
Trichur – 680 009, India.
*sebastian_sherly@yahoo.com*

**Abstract.** Here we develop a first order autoregressive model $\{X_n\}$ that is marginally stationary where $X_n$ is the sum/ extreme of $k$ i.i.d observations. We prove that stationary solutions to these models are also either semi-selfdecomposable/ extreme-semi-selfdecomposable or, sum/ extreme stable with respect to Harris distribution.

## 1. Introduction

A sequence $\{X_n\}$ of random variables (r.v) describes the additive first order autoregressive (AR(1)) scheme considered here if there exists an innovation sequence $\{\varepsilon_n\}$ of independent and identically distributed (i.i.d) r.vs satisfying

$$X_n = bX_{n-1} + \varepsilon_n , \ \forall \ n > 0 \text{ integer and some } 0 < b < 1. \tag{1.1}$$

$\{X_n\}$ is marginally stationary if $X_n \overset{d}{=} X_{n-1} \ \forall \ n > 0$ integer.

This investigation is motivated by the possibility whereby the AR(1) sequence $\{X_n\}$ is composed of $k$ independent AR(1) sequences $\{Y_{i,n}\}$ , $i = 1, 2, \ldots., k$, and where for each $n > 0$ integer $Y_{i,n}$ , $i = 1, 2, \ldots., k$ are identically distributed. That is, for each $n$, $X_n = \sum_{i=1}^{k} Y_{i,n}$ and $Y_{i,n}$, $i = 1,2, \ldots, k$ are i.i.d, $k$ being a fixed positive integer. For example, the variable $X_n$ could be the quantity of water flowing through a river, or the number of patients in a hospital, or the sales of a particular item by an agency, or the number of items produced in a factory having more than one plant for the production of the same. In all these cases the resultant observation $X_n$ is either the sum of the quantities $Y_{i,n}$, $i = 1,2, \ldots, k$ of water flowing through $k$ tributaries of the river, or the sum of the number of patients $Y_{i,n}$, $i = 1,2, \ldots, k$ in $k$ different specialities in the hospital, or the sum of the sales $Y_{i,n}$, $i = 1,2, \ldots, k$ by the agency through their $k$ different retail outlets or the sum of the quantities produced $Y_{i,n}$, $i = 1,2, \ldots, k$ at the $k$ different plants in the factory. Thus we generalize the AR(1) model (1.1) where the observations $X_n$ are the sum of $k$ i.i.d r.vs, and then extend the discussion to the maximum and minimum schemes. We need the following concepts in our discussion.

**Definition.1.1** (Pillai, 1985). A characteristic function (CF) $f$ is semi-$\alpha$-Laplace($a,b$) if $\forall$ $t \in \mathbf{R}$ and for some $0 < b < 1 < a$,

$$f(t) = 1/\{1 + \psi(t)\}, \text{ where } \psi(t) = a\psi(bt), ab^\alpha = 1, \alpha \in (0,2].$$

Satheesh, *et al.* (2002) has considered a generalized semi-$\alpha$-Laplace($a,b,k$) law having CF $\{1 + \psi(t)\}^{-1/k}$ where $\psi(t)$ is as above. The generalized semi Mittag-Leffler($a,b,k$) (ML) laws with Laplace transform (LT) $\varphi(s) = \{1 + \psi(s)\}^{-1/k}$, $s > 0$, $\psi(s) = a\psi(bs)$, $0 < b < 1 < a$, $ab^\alpha = 1$, $\alpha \in (0,1]$ is its non-negative analogue and the generalized discrete semi-ML($a,b,k$) laws having probability generating function (p.g.f) $P(s) = \{1 + \psi(1-s)\}^{-1/k}$, $0 < s < 1$, where $\psi(1-s) = a\psi(b(1-s))$, $0 < b < 1 < a$, $ab^\alpha = 1$, $\alpha \in (0,1]$ is its non-negative integer-valued analogue, see Satheesh, *et al.* (2002). In these families only the generalized Laplace law with CF $\{1 + \lambda t^2\}^{-1/k}$ has finite variance. Similarly, in the non-negative case only the gamma($1/k, \lambda$) law has finite mean and in the discrete case only the negative binomial($1/k, \lambda$) law has finite mean. Also if $\psi(t) = a\psi(bt)$ for two values of $b$, say $b_1$ and $b_2$ such that $ln(b_1)/ln(b_2)$ is irrational then $\psi(t) = \lambda|t|^\alpha$, $\lambda > 0$.

**Definition.1.2** (Maejima and Naito, 1998). A CF $f$ is semi-selfdecomposable($b$) (SSD($b$)) if for some $0 < b < 1$ there exists a CF $f_o$ that is infinitely divisible (ID) such that

$$f(t) = f(bt) f_o(t), \forall \ t \in \mathbf{R}.$$

If this relation holds for every $0 < b < 1$ then $f$ is selfdecomposable (SD).

**Definition.1.3** (Megyesi, 2002). A distribution function (d.f) $F$ is max-semi-stable($a,c$) if either,

$$F(x) = exp\{-x^{-\alpha} h(ln(x))\}, x > 0, \alpha > 0,$$

where $h(x)$ is a positive bounded periodic function with period $ln(c)$, $c>1$, and there exists an $a>1$ such that $ac^{-\alpha} = 1$, or,

$$F(x) = exp\{-|x|^{\alpha} h(ln(|x|))\}, x<0, \alpha>0,$$

where $h(x)$ is as above with period $|ln(c)|$, $c<1$, and there is an $a>1$ such that $ac^{\alpha} = 1$.

These are the extended Frechet and the extended Weibull types. The extended Gumbel has the translation invariance property and is not considered here.

**Remark.1.1** The first d.f in definition.1.3 can be written in the form $exp\{-\psi(x)\}$, where $\psi(x)$ satisfies $\psi(x) = a\psi(cx)$, $x>0$ for some $a>1$, $c>1$, and there is an $\alpha>0$ satisfying $ac^{-\alpha} =1$. Similarly, the second one also, where $\psi(x) = a\psi(cx)$, $x<0$, for some $a>1$, $c<1$, and $\alpha>0$ satisfying $ac^{\alpha} = 1$, see Satheesh and Sandhya (2006$b$).

**Definition.1.4** (Becker-Kern, 2001). A non-degenerate d.f $F$ is max-SSD($c$) if for some $c>1$ and $v\in\mathbf{R}$ there is a non-degenerate d.f $H$ such that

$$F(x) = F(c^{v}x + \beta) H(x), \forall x\in\mathbf{R},$$

where $\beta=0$ if $v\neq0$ and $\beta = ln(c)$ if $v=0$. Here we will consider the case of $\beta=0$, ie. $v\neq0$ only so that the above relation becomes $F(x) = F(cx) H(x), \forall x\in\mathbf{R}$, and some $c\in(0,1)\cup(1,\infty)$. If the relation holds for every $c\in(0,1)\cup(1,\infty)$, then $F$ is max-SD.

Satheesh and Sandhya (2005) showed that the sequence $\{X_n\}$ of r.vs generates the marginally stationary AR(1) scheme (1.1) if $X_n$ is SSD($b$) and also discussed the integer-valued analogue of it. Extending the discussion to the maximum scheme Satheesh and Sandhya (2006$b$) showed that $\{X_n\}$ generates the marginally stationary max-AR(1) scheme iff $X_n$ is max-SSD($c$). The structure here is:

$$X_n = bX_{n-1}\vee \varepsilon_n , \forall n>0 \text{ integer and some } b>0, \text{ and } c=1/b. \tag{1.2}$$

They then modified the max-AR(1) scheme to:

$$\left. \begin{array}{l} X_n \ = bX_{n-1} , \text{ with probability } p \\ \ \ \ \ = bX_{n-1}\vee \varepsilon_n , \text{ with probability } (1\text{-}p). \end{array} \right\} \tag{1.3}$$

Subsequently Satheesh and Sandhya (2006$b$) showed, assuming $X_0 \overset{d}{=} \varepsilon_1$, that a sequence $\{X_n\}$ of r.vs generates the marginally stationary max-AR(1) scheme (1.3), iff $X_n$ is exponential max-semi-stable($a,c$), $a=\frac{1}{p}$, $c=\frac{1}{b}$ having d.f of the form:

$$F(x) = \frac{1}{1+\psi(x)} = \frac{1}{1+\frac{1}{p}\psi(cx)} , \forall x\in\mathbf{R},$$

where $\psi(x)$ is as in remark.1.1 with $c{>}1$ if $X_n{>}0$ and $c{<}1$ if $X_n{<}0$ $\forall n$.

In the additive AR(1) and minimum AR(1) schemes analogous to (1.3), Jayakumar and Pillai (1992), Pillai (1991), Balakrishna and Jayakumar (1997) have similar results.

Whether additive, minimum or maximum; the structure in (1.3) essentially captures a geometric sum or a geometric extreme scheme, the geometric law being supported by the set $\{1, 2, 3, \ldots.\}$. A generalization of the geometric law is the Harris(1,$a$,$k$) law on $\{1, 1{+}k, 1{+}2k, \ldots.\}$, that is described by its p.g.f

$$P(s) = \frac{s}{\{a - (a-1)s^k\}^{1/k}} , k{>}0 \text{ integer and } a{>}1.$$

The stability properties of the Harris(1,$a$,$k$) law in the summation scheme were studied by Satheesh, *et al*. (2002) while Satheesh and Nair (2002$b$, 2004) studied them for the minimum and maximum. More distributional and divisibility properties, simulation and estimation problems of Harris(1,$a$,$k$) law have been addressed to in Sandhya, *et al*. (2006). As it turns out, the stationary solution to the generalized AR(1) model that we discuss here has the property of Harris-sum/ extreme stability.

Thus in the next section we consider a generalization of (1.1) where for each $n{>}0$, $X_n \stackrel{d}{=} \sum_{i=1}^{k} Y_{i,n}$ and $Y_{i,n}$, $i =1,2, \ldots., k$ are i.i.d, $k$ being a fixed positive integer. The case of the maximum scheme instead of addition is also considered along with this. The discussion is then extended to the addition scheme for integer-valued $X_n$ in section.3. An off-shoot of the development here is the description of SSD laws on $\{0, m, 2m, \ldots.\}$, $m{>}0$ integer. In section.4 we discuss the case for the minimum scheme. In section.5 we give a different formulation of the main results.

## 2. Generalizations of the Additive and Maximum AR(1) Schemes.

We now modify the model (1.1) as follows. We have $k$ independent AR(1) sequences $\{Y_{i,n}\}$, $i = 1,2, \ldots, k$, where $Y_{i,n} = bY_{i,n-1} + \varepsilon_{i,n}$, $\forall$ $n{>}0$ integer and some $0{<}b{<}1$, where $Y_{i,n}$, $i = 1,2, \ldots., k$ are identically distributed. Here, $\varepsilon_{i,n}$, $i = 1,2, \ldots., k$, $n{>}0$ integer are i.i.d. Hence (1.1) is described in terms of $\{Y_{i,n}\}$ as;

$$\sum_{i=1}^{k} Y_{i,n} = b\sum_{i=1}^{k} Y_{i,n-1} + \sum_{i=1}^{k} \varepsilon_{i,n} , \forall n{>}0 \text{ integer and some } 0{<}b{<}1. \qquad (2.1)$$

Assuming $\{Y_{i,n}\}$ to be marginally stationary (that is, $Y_{i,n} \stackrel{d}{=} Y_{i,n-1}$ $\forall$ $n$) their CFs satisfy;

$$f_y^{\,k}(u) = f_y^{\,k}(bu)\ f_\varepsilon^{\,k}(u)\text{, for some } 0<b<1,$$

$$= \{\, f_y(bu)\ f_\varepsilon(u)\,\}^k\text{. Hence;}$$

**Theorem.2.1** A sequence $\{Y_{i,n}\}$ of r.vs defines the marginally stationary AR(1) scheme (2.1) if $Y_{i,n}$ is SSD($b$).

**Corollary.2.1** A sequence $\{Y_{i,n}\}$ describes the marginally stationary AR(1) scheme (2.1) for all $0<b<1$ if $Y_{i,n}$ is SD.

In the maximum scheme the model (2.1) reads as:

$$\bigvee_{i=1}^{k} Y_{i,n} = b\{ \bigvee_{i=1}^{k} Y_{i,n-1} \} \vee \{ \bigvee_{i=1}^{k} \varepsilon_{i,n} \}\text{, } \forall\ n>0 \text{ integer and some } b>0. \qquad (2.2)$$

Assuming $\{Y_{i,n}\}$ to be marginally stationary, in terms of their d.fs this equation reads:

$$F^k(x) = F^k(cx)\ G^k(x)\text{, } c=1/b\text{. Or}$$

$$= \{\, F(cx)\ G(x)\,\}^k\text{, for some } c>0\text{. Hence;}$$

**Theorem.2.2** A sequence $\{Y_{i,n}\}$ of r.vs defines the marginally stationary max-AR(1) scheme (2.2) if $Y_{i,n}$ is max-SSD($c$), $c=1/b$. $\{Y_{i,n}\}$ describes (2.2) $\forall\ b>0$ if $Y_{i,n}$ is max-SD.

Now modifying (2.1) we may have, $\forall\ n>0$ integer and some $0<b<1$:

$$\sum_{i=1}^{k} Y_{i,n} = b\sum_{i=1}^{k} Y_{i,n-1} \qquad \text{with probability } p$$

$$\sum_{i=1}^{k} Y_{i,n} = b\sum_{i=1}^{k} Y_{i,n-1} + \sum_{i=1}^{k} \varepsilon_{i,n} \qquad \text{with probability } (1\text{-}p) \qquad (2.3)$$

Assuming $Y_{i,o} \overset{d}{=} \varepsilon_{i,1}$ and marginal stationarity of $\{Y_{i,n}\}$ their CFs satisfy; for $n=1$:

$\phi^k(t) = p\phi^k(bt) + (1\text{-}p)\ \phi^k(bt)\ \phi^k(t)$. That is;

$\phi^k(t) = p\phi^k(bt)/\{1\text{-}(1\text{-}p)\phi^k(bt)\}$. Or;

$\phi^k(t) = \phi^k(bt)/\{a\text{-}(a\text{-}1)\phi^k(bt)\}$, where $a=1/p$. Hence;

$\phi(t) = \{\phi^k(bt)/\{a\text{-}(a\text{-}1)\phi^k(bt)\}\}^{1/k}$. $\qquad (*)$

Equation $(*)$ means that $Y_{i,1}$ is Harris(1,$a$,$k$)-sum stable. Hence by the characterization of generalized semi-$\alpha$-Laplace laws, Satheesh, *et al.* (2002, theorem.2.1), $Y_{i,1}$ is generalized semi-$\alpha$-Laplace($a$,$b$,$k$) with CF $\phi(t) = \{1+\psi(t)\}^{-1/k}$ where $\psi(t)$ satisfies $\psi(t) = a\,\psi(bt)$, $a$ and $k$ being that in Harris(1,$a$,$k$) and $ab^\alpha = 1$ for some $\alpha \in (0,2]$, $a=1/p$ in the model (2.3). That is, $Y_{i,1}$, $Y_{i,o}$ and $\varepsilon_{i,1}$ are generalized semi-$\alpha$-Laplace($a$,$b$,$k$). Now, by the marginal stationarity of $\{Y_{i,n}\}$ we get $Y_{i,2}$ is generalized semi-$\alpha$-Laplace($a$,$b$,$k$) and recursively $Y_{i,3}$, $Y_{i,4}$, …. also follow generalized semi-$\alpha$-Laplace($a$,$b$,$k$). Such an induction argument over the index $n$ results in:

**Theorem.2.3** Under the assumption $Y_{i,o} \overset{d}{=} \varepsilon_{i,1}$, a sequence $\{Y_{i,n}\}$ of r.vs defines the marginally stationary AR(1) scheme (2.3) iff $Y_{i,n}$ is generalized semi-$\alpha$-Laplace($\frac{1}{p}$,$b$,$k$).

The max-analogue of the characterization of generalized semi-$\alpha$-Laplace laws mentioned above is the following. The proof follows on the same lines as that of theorem.2.1 in Satheesh, *et al.* (2002).

**Theorem.2.4** The i.i.d r.vs $X_i$ are Harris(1,$a$,$k$)-max stable iff its d.f is $F(x) = \{1+\psi(x)\}^{-1/k}$, where $\psi(x)$ satisfies $\psi(x) = a\,\psi(cx)$, $a$ and $k$ being that in Harris(1,$a$,$k$) and $ac^\alpha = 1$ for some $c>0$ and $\alpha>0$.

Satheesh and Sandhya (2006$a$) had discussed $\varphi$-max-semi-stable laws for a LT $\varphi$. In this terminology the above d.f is called a gamma-max-semi-stable($a$,$c$,1/$k$) law. When $k=1$ we have the exponential max-semi-stable($a$,$c$) model that characterized (1.3).

Now, modifying (2.2) further we have, $\forall$ $n>0$ integer and some $b>0$:

$$
\left.
\begin{aligned}
&\overset{k}{\underset{i=1}{\bigvee}} Y_{i,n} = b\{ \overset{k}{\underset{i=1}{\bigvee}} Y_{i,n-1}\} && \text{with probability } p \\[2ex]
&\overset{k}{\underset{i=1}{\bigvee}} Y_{i,n} = b\{ \overset{k}{\underset{i=1}{\bigvee}} Y_{i,n-1}\} \vee \{ \overset{k}{\underset{i=1}{\bigvee}} \varepsilon_{i,n}\} && \text{with probability } (1\text{-}p)
\end{aligned}
\right\} \qquad (2.4)
$$

Proceeding as in the additive scheme we have the max-analogue of theorem.2.3.

**Theorem.2.5** Under the assumption $Y_{i,o} \overset{d}{=} \varepsilon_{i,1}$, a sequence $\{Y_{i,n}\}$ of r.vs defines the marginally stationary max-AR(1) scheme (2.4) iff $Y_{i,n}$ is gamma-max-semi-stable($\frac{1}{p}$,$\frac{1}{b}$,$\frac{1}{k}$).

Suppose we require (2.3) to be satisfied for all $b \in (0,1)$, then:

**Corollary.2.2** Assuming $Y_{i,o} \overset{d}{=} \varepsilon_{i,1}$, a sequence $\{Y_{i,n}\}$ of r.vs defines the marginally stationary AR(1) scheme (2.3) for all $b \in (0,1)$, iff $Y_{i,n}$ is generalized Linnik (generalized $\alpha$-Laplace) with CF $\{1+\lambda|t|^\alpha\}^{-1/k}$, $k>0$ integer, $\alpha \in (0,2]$ and $\lambda>0$.

Suppose we further demand $Y_{i,n}$ to have finite variance, then (2.3) characterizes the generalized Laplace law with CF $\{1+\lambda t^2\}^{-1/k}$, $k>0$ integer and $\lambda>0$.

Suppose we need $Y_{i,n}>0$, then by the Harris-sum stability characterization of generalized semi-ML laws, Satheesh, *et al.* (2002, corollary.2.1), we have:

**Theorem.2.6** Under the assumption $Y_{i,o} \overset{d}{=} \varepsilon_{i,1}$, a sequence $\{Y_{i,n}>0\}$ of r.vs defines the marginally stationary AR(1) scheme (2.3), iff $Y_{i,n}$ is generalized semi-ML($\frac{1}{p}$, $b$, $k$).

Under the additional assumptions as in corollary.2.2 we have:

**Corollary.2.3** Assuming $Y_{i,o} \overset{d}{=} \varepsilon_{i,1}$, a sequence $\{Y_{i,n}>0\}$ of r.vs defines the marginally stationary AR(1) scheme (2.3) for all $b \in (0,1)$, iff $Y_{i,n}$ is generalized ML with LT $\{1+\lambda s^\alpha\}^{-1/k}$, $k>0$ integer, $\alpha \in (0,1]$ and $\lambda>0$. Further demanding $Y_{i,n}$ to have finite mean the gamma($1/k,\lambda$) law with LT $\{1+\lambda s\}^{-1/k}$, $k>0$ integer and $\lambda>0$ is characterized by (2.3).

### 3. Generalization of the Discrete AR(1) Scheme and SSD Laws with Gaps

To discuss (2.3) for discrete r.vs we need the description of integer-valued r.vs of the same type in Satheesh and Nair (2002*a*) that we present here as a remark.

**Remark.3.1** If $\phi(s)$ is a LT, then $P(s) = \phi(1-s)$, $0<s<1$ is a p.g.f. If $\phi_1$ and $\phi_2$ are LTs, then the p.g.fs $P_1(s) = \phi_1(1-s)$ and $P_2(s) = \phi_2(1-s)$ are of the same type iff $\phi_1(1-s) = \phi_2(c(1-s))$, for all $0<s<1$ and some $c>0$. Two p.g.fs $P_1(s)$ and $P_2(s)$ are of the same type iff $P_1$ is a $P_2$ compounded Bernoulli law. Thus in the setup of integer-valued r.vs the equivalent of r.vs of the same type is obtained by replacing $bX$ by $b \circ X = \sum_{i=1}^{X} Z_i$, where $\{Z_i\}$ are i.i.d Bernoulli($b$) r.vs independent of $X$ with $P\{Z_i=0\} = 1-b$.

Consequently (2.3) becomes:

$$\sum_{i=1}^{k} Y_{i,n} = b \circ \sum_{i=1}^{k} Y_{i,n-1} \qquad \text{with probability } p$$

$$\sum_{i=1}^{k} Y_{i,n} = b \circ \sum_{i=1}^{k} Y_{i,n-1} + \sum_{i=1}^{k} \varepsilon_{i,n} \quad \text{with probability } (1-p) \tag{3.1}$$

Discrete generalized semi-ML($a,b,k$) laws with p.g.f $\{1+\psi(1-s)\}^{-1/k}$, $0<s<1$, where $\psi(1-s)$ satisfies $\psi(1-s) = a\psi(b(1-s))$, $ab^\alpha = 1$ for some $\alpha \in (0,1]$, were characterized in theorem.3.3 of Satheesh, *et al.* (2002) by the property of Harris-sum stability of distributions on {0, 1, 2, ….}. Hence when {$Y_{i,n}$} are integer-valued we have:

**Theorem.3.1** Under the assumption $Y_{i,o} \overset{d}{=} \varepsilon_{i,1}$, a sequence {$Y_{i,n}$} of non-negative integer-valued r.vs defines the marginally stationary AR(1) scheme (3.1), iff $Y_{i,n}$ is discrete generalized semi-ML($\frac{1}{p}$, $b$, $k$).

Now we characterize discrete generalized ML and negative binomial (NB) laws respectively by additional conditions as in corollary.2.2 on the scheme.

**Corollary.3.1** Assuming $Y_{i,o} \overset{d}{=} \varepsilon_{i,1}$, a sequence {$Y_{i,n}$} of non-negative integer-valued r.vs defines the marginally stationary AR(1) scheme (3.1) for all $b \in (0,1)$, iff $Y_{i,n}$ is discrete generalized ML with p.g.f $\{1+\lambda(1-s)^{\alpha}\}^{-1/k}$, $k>0$ integer, $\alpha \in (0,1]$ and $\lambda>0$. Further demanding $Y_{i,n}$ to have finite mean the model (3.1) characterizes the NB(1/$k$, $\lambda$) law with p.g.f $\{1+\lambda(1-s)\}^{-1/k}$, $k>0$ integer and $\lambda>0$.

We next characterize distributions that have gaps in its support by Harris-sum stability as follows. Certain implications of having and not having gaps in the support were investigated in Satheesh (2004).

**Theorem.3.2** A distribution on {0, $m$, 2$m$, 3$m$, ….}, $m>0$ integer, is Harris(1, $a$, $k$)-sum stable iff its p.g.f is $P(s) = \{1+\psi(1-s^m)\}^{-1/k}$, $\psi(1-s^m)$ satisfying $\psi(1-s^m) = a\psi(b(1-s^m))$, $ab^{\alpha} =1$ for some $\alpha \in (0,1]$.

Proof. The assertion follows from lemma.4.1 and theorem.3.3 in Satheesh, *et al.* (2002).

Now in terms of AR(1) models we have:

**Theorem.3.3** Under the assumption $Y_{i,o} \overset{d}{=} \varepsilon_{i,1}$, a sequence {$Y_{i,n}$} of r.vs on {0, $m$, 2$m$, 3$m$, ….}, $m>0$ integer, defines the marginally stationary AR(1) scheme (3.1), iff the p.g.f of $Y_{i,n}$ is $P(s) = \{1+\psi(1-s^m)\}^{-1/k}$, $k>0$ integer, $p\psi(1-s^m) = \psi(b(1-s^m))$, $b^{\alpha} = p$, $\alpha \in (0,1]$, $0<s<1$.

Under the additional conditions on the scheme as in corollary.2.2 we have:

**Corollary.3.2** Assuming $Y_{i,o} \overset{d}{=} \varepsilon_{i,1}$, a sequence {$Y_{i,n}$} of r.vs on {0, $m$, 2$m$, 3$m$, ….}, $m>0$ integer, defines the marginally stationary AR(1) scheme (3.1) $\forall b \in (0,1)$, iff the p.g.f of $Y_{i,n}$ is $\{1+\lambda(1-s^m)^{\alpha}\}^{-1/k}$, $k>0$ integer, $\alpha \in (0,1]$ and $\lambda>0$. Further demanding $Y_{i,n}$ to have finite mean the model (3.1) characterizes the p.g.f of $Y_{i,n}$ as $\{1+\lambda(1-s^m)\}^{-1/k}$, $k>0$ integer and $\lambda>0$.

**Remark.3.2** Recently Satheesh and Sandhya (2005) have described SSD p.g.fs. Incidentally the above discussion suggests possible extension of the notion of SD and SSD laws to distributions on $\{0, m, 2m, 3m, \ldots\}$, $m>0$ integer.

**Definition.3.1** An integer-valued distribution on $\{0, m, 2m, 3m, \ldots\}$, $m>0$ integer, with p.g.f $P(s^m)$ is discrete SSD($b$) if for some $0<b<1$, there exists another p.g.f $P_o(s^m)$ that is ID such that

$$P(s^m) = P(1-b+bs^m) \, P_o(s^m), \, \forall \, s \in (0,1).$$

The distribution is SD if the above relation holds good for all $0<b<1$. If the p.g.fs are described in terms of LTs as in remark.3.1 above, then equivalently we have: A p.g.f $P(s^m) = \phi(1-s^m)$ is discrete SSD($b$) if for some $0<b<1$, there exists another p.g.f $P_o(s^m) = \phi_o(1-s^m)$ that is ID such that

$$\phi(1-s^m) = \phi[b(1-s^m)] \, \phi_o(1-s^m).$$

**Example.3.1** By the Harris(1,$a$,$k$)-sum stability of the p.g.f $\{1+\psi(1-s^m)\}^{-1/k}$ (see theorem.3.2) we get;

$$\frac{1}{\{1+\psi(1-s^m)\}^{1/k}} = \frac{1}{\{1+\psi[b(1-s^m)]\}^{1/k}} \frac{1}{\{a-(a-1)/\{1+\psi[b(1-s^m)]\}\}^{1/k}}.$$

Here the second factor itself is a p.g.f being the Harris-sum of the distribution with p.g.f $\{1+\psi[b(1-s^m)]\}^{-1/k}$, this Harris law being supported on $\{0, k, 2k, \ldots\}$. This Harris-sum is also ID since the Harris law is ID. Hence the p.g.f $\{1+\psi(1-s^m)\}^{-1/k}$ is discrete SSD($b$) on $\{0, m, 2m, 3m, \ldots\}$, $m>0$ integer.

**Theorem.3.4** A sequence $\{X_n\}$ of integer-valued r.vs defines a marginally stationary AR(1) sequence on $\{0, m, 2m, 3m, \ldots\}$, $m>0$ integer, with $0<b<1$ if $X_n$ is discrete SSD($b$) on $\{0, m, 2m, 3m, \ldots\}$, $m>0$ integer.

Proof. The proof follows on lines similar to that of theorem.1 in Satheesh and Sandhya (2005).

## 4. Generalizations of the Minimum AR(1) Scheme.

To extend the model to the minimum scheme we need describe min-SSD($c$) laws. Here we do not attempt a detailed study of this class and also restrict our discussion to the support $[0,\infty)$. More on this class is in Satheesh and Sandhya (2006$c$).

**Definition.4.1** A non-degenerate d.f $F$ with survival function (s.f) $R$ is min-SSD($c$) if for some $0<c<1$ there is another s.f $S$ such that

$$R(x) = R(cx) \, S(x), \, \forall \, x>0.$$

If this is true for every $0<c<1$, then $F$ is min-SD.

The following results also give certain examples in these classes.

**Theorem.4.1** Generalized semi-Pareto$(p,\alpha,1/k)$ law with s.f $\{1+\psi(x)\}^{-1/k}$, $x>0$, $p\,\psi(x) = \psi(p^{1/\alpha}x)$, $\forall x>0$, some $0<p<1$, and $\alpha>0$ is min-SSD$(p^{1/\alpha})$.

Proof. By the Harris$(1,a,k)$-min stability of the generalized semi-Pareto$(p,\alpha,1/k)$ law, (Satheesh and Nair (2002$b$), quoted below as theorem.4.3) it follows that we can write:

$$\{1+\psi(x)\}^{-1/k} = \{1+\psi(bx)\}^{-1/k}/\{a-(a-1)\,[1+\psi(bx)]^{-1}\}^{1/k}, \; p = 1/a, \; b{=}p^{1/\alpha}.$$

Now the assertion follows as done in example.3.1.

A more general approach to proving that a distribution is SSD or max/ min SSD is from the angle of mixtures as done in Satheesh and Sandhya (2005, 2006$b$, $c$).

**Corollary.4.1** Semi-Pareto$(p,\alpha)$ laws of Pillai (1991) are min-SSD$(p^{1/\alpha})$, which follows by its geometric-min stability, the geometric law being on $\{1, 2, \ldots.\}$.

**Corollary.4.2** Pareto laws with s.f $1/\{1+s^{\alpha}\}$, $\alpha>0$, are min-SD.

Now we modify the scheme (2.2) to the minimum structure as:

$$\bigwedge_{i=1}^{k} Y_{i,n} = b\{\bigwedge_{i=1}^{k} Y_{i,n-1}\} \wedge \{\bigwedge_{i=1}^{k} \varepsilon_{i,n}\} \; \forall \, n>0 \text{ integer and some } b>1. \qquad (4.1)$$

Assuming $\{Y_{i,n}\}$ to be marginally stationary, in terms of the s.fs $R$ of $Y_{i,n}$ and $S$ of $\varepsilon_{i,n}$ this equation reads:

$$R^k(x) = R^k(cx)\,S^k(x)\,, c{=}1/b. \text{ Or}$$

$$= \{R(cx)\,S(x)\}^k, \text{ for some } 0<c<1. \text{ Hence;}$$

**Theorem.4.2** A sequence $\{Y_{i,n}\}$ of non-negative r.vs defines a marginally stationary min-AR(1) scheme as in (4.1) if $Y_{i,n}$ is min-SSD$(c)$, $c{=}1/b$. $\{Y_{i,n}\}$ describes the structure (4.1) for all $b>1$ if $Y_{i,n}$ is min-SD.

The following is a restatement of proposition.3 in Satheesh and Nair (2002$b$).

**Theorem.4.3** The i.i.d r.vs $X_i$ are Harris$(1,a,k)$-min stable iff it is generalized semi-Pareto$(p,\alpha,1/k)$ with d.f $F(x) = 1-\{1+\psi(x)\}^{-1/k}$, $x>0$, $p\,\psi(x) = \psi(p^{1/\alpha}x)$, $\forall\, x>0$, some $0<p<1$, $p = 1/a$, and $\alpha>0$.

Now, modifying (4.1) further we have, $\forall n>0$ integer and some $b>1$:

$$\bigwedge_{i=1}^{k} Y_{i,n} = b\{ \bigwedge_{i=1}^{k} Y_{i,n-1} \} \qquad \text{with probability } p$$

$$\bigwedge_{i=1}^{k} Y_{i,n} = b\{ \bigwedge_{i=1}^{k} Y_{i,n-1} \} \wedge \{ \bigwedge_{i=1}^{k} \varepsilon_{i,n} \} \text{ with probability } (1\text{-}p) \qquad (4.2)$$

Assuming marginal stationarity of $Y_{i,n}$ and $Y_{i,o} \overset{d}{=} \varepsilon_{i,1}$, and proceeding as in the additive scheme we have the following min-analogue of theorem.2.3.

**Theorem.4.4** Under the assumption $Y_{i,o} \overset{d}{=} \varepsilon_{i,1}$, a sequence $\{Y_{i,n}\}$ of non-negative r.vs defines the marginally stationary min-AR(1) scheme (4.2), iff $Y_{i,n}$ is generalized semi-Pareto$(p,\alpha,1/k)$ and $b= p^{-1/\alpha}$.

## 5. Concluding Remarks.

Notice that the schemes (2.2) and (2.4) include $b{<}1$ and the explosive case $b{>}1$ as well. This is also in tune with the range of the parameter $c$ in max-SSD and max-semi-stable laws. Finally, a different formulation of theorem.2.3 in terms of the innovation sequence is:

**Theorem.5.1** Under the assumption $Y_{i,o} \overset{d}{=} \varepsilon_{i,1}$, the innovation sequence $\{\varepsilon_{i,n}\}$ describes the marginally stationary AR(1) scheme $\{Y_{i,n}\}$ in (2.3) iff $\varepsilon_{i,1}$ is generalized semi-$\alpha$-Laplace$(\frac{1}{p},b,k)$.

Proof. The only additional thing to be proved here from theorem.2.3 is the marginal stationarity of the scheme. Under the given assumptions the CF $f_1$ of $Y_{i,1}$ is given by:

$$f_1^{k}(t) = \frac{p}{1+\psi(bt)} + \frac{1-p}{[1+\psi(bt)][1+\psi(t)]} = \frac{1}{1+\psi(t)} \text{, since } p\,\psi(t) = \psi(bt).$$

Hence by induction the scheme (2.3) is marginally stationary, completing the proof.

Similar formulations of theorems 2.5 and 4.4 in the maximum and minimum schemes in terms of their innovation sequence $\{\varepsilon_{i,n}\}$ are as follows. Proofs follow on similar lines.

**Theorem.5.2** Under the assumption $Y_{i,o} \overset{d}{=} \varepsilon_{i,1}$, the innovation sequence $\{\varepsilon_{i,n}\}$ describes the marginally stationary max-AR(1) scheme $\{Y_{i,n}\}$ in (2.4) iff $\varepsilon_{i,1}$ is gamma-max-semi-stable$(\frac{1}{p}, \frac{1}{b}, \frac{1}{k})$.

**Theorem.5.3** Under the assumption $Y_{i,o} \stackrel{d}{=} \varepsilon_{i,1}$ , the innovation sequence $\{\varepsilon_{i,n}\}$ describes the marginally stationary min-AR(1) scheme $\{Y_{i,n}\}$ in (4.2) iff $\varepsilon_{i,1}$ is generalized semi-Pareto$(p,\alpha,1/k)$.

Similar formulations of theorems 3.1 and 3.3 result in characterizations of discrete generalized semi-ML law and that with p.g.f $\{1+\psi(1-s^m)\}^{-1/k}$. Thus we have potentially useful generalized AR(1) schemes in the additive, maximum and minimum structures.

**Acknowledgements:** Authors thank the referee whose careful reading corrected certain ambiguous statements and typographical errors and also the executive editor for a better presentation of the material. The third author's work is supported by a fellowship from the UGC, India.